\documentclass[11pt,a4paper]{article}
\usepackage{a4wide}
\usepackage{latexsym}
\usepackage{amsmath}
\usepackage{amsthm}
\usepackage{amssymb,enumerate}
\theoremstyle{plain}
\newtheorem{theo}{Theorem}[section]
\newtheorem{lemma}[theo]{Lemma}
\newtheorem{prop}[theo]{Proposition}
\newtheorem{coro}[theo]{Corollary}
\theoremstyle{definition}
\newtheorem{defi}[theo]{Definition}
\newtheorem{rema}[theo]{Remark}
\newtheorem{exam}[theo]{Example}

\newenvironment{proof1}{\medskip\par\noindent{\bf Proof.}}{\hfill $\Box$
\medskip\par}

\def\C{\mathbb{C}}
\def\N{\mathbb{N}}
\def\R{\mathbb{R}}

\def\a{\alpha}
\def\b{\beta}
\def\ro{\rho}
\def\ga{\gamma}

\def\bmu{\boldsymbol{\mu}}
\def\bM{\mathbb{M}}
\def\M{\mathbb{M}}
\def\bm{\boldsymbol{m}}

\def\parn{\par\noindent}

\begin{document}

\title{Flat functions in Carleman ultraholomorphic classes\\via proximate orders}
\author{Javier Sanz}
\date{\today}

\maketitle

{ \small
\begin{center}
{\bf Abstract}
\end{center}

Whenever the defining sequence of a Carleman ultraholomorphic class (in the sense of H. Komatsu) is strongly regular and associated with a proximate order,
flat functions are constructed in the class on sectors of optimal opening. As consequences, we obtain analogues of both Borel--Ritt--Gevrey theorem and Watson's lemma in this general situation.\par
}

\bigskip \bigskip

\section{Introduction}

The Carleman ultraholomorphic classes $\mathcal{A}_{\M}(S)$ in a sector $S$ of the Riemann surface of the logarithm consist of those holomorphic functions $f$ in $S$ whose derivatives of order $n\ge 0$ are uniformly bounded there by, essentially, the values $n!M_n$, where $\M=(M_n)_{n\in\N_0}$ is a sequence of positive real numbers. In case bounds are not uniform on $S$ but are valid and depend on every proper subsector of $S$ to which the function is restricted, we obtain the class $\tilde{\mathcal{A}}_{\M}(S)$ of functions with a (non-uniform) $\M$-asymptotic expansion at 0 in $S$, given by a formal power series $\hat f=\sum_{n\ge 0}a_nz^n/n!$ whose coefficients are again suitably bounded in terms of $\M$ (we write $f\sim_{\M}\hat f$ and $(a_n)_{n\in\N_0}\in\Lambda_{\M}$). The map sending $f$ to $(a_n)_{n\in\N_0}$ is the asymptotic Borel map $\tilde{\mathcal{B}}$, and  $f$ is said to be flat if $\tilde{\mathcal{B}}(f)$ is the null sequence. See Subsection~\ref{subsectCarlemanclasses} for the precise definitions of all these classes and concepts.\par

In order to obtain good properties for these classes, the sequence $\M$ is usually subject to some standard conditions; in particular, we will only consider strongly regular sequences as defined by V. Thilliez~\cite{thilliez}, see Subsection~\ref{subsectstrregseq}.
The best known example is that of Gevrey classes, appearing when the sequence is chosen to be $\bM_{\a}=(n!^{\a})_{n\in\N_{0}}$, $\a>0$, and for which we use the notations $\mathcal{A}_{\a}(S)$, $\tilde{\mathcal{A}}_{\a}(S)$, $\Lambda_{\a}$, $f\sim_{\a}\hat{f}$ and so on, for simplicity. Let us denote by $S_\ga$ the sector bisected by the direction $d=0$ and with opening $\pi\ga$. It is well known that $\tilde{\mathcal{B}}:\tilde{\mathcal{A}}_{\a}(S_\ga)\to\Lambda_{\a}$ is surjective if, and only if, $\ga\le\a$ (Borel--Ritt--Gevrey theorem, see~\cite{Ramis1,Ramis2,MartinetRamis}, \cite[Thm.\ 2.2.1]{Balser}). It is natural to call this an extension result, and to think then about the possibility of obtaining linear and continuous right inverses for $\tilde{\mathcal{B}}$ in suitably topologized classes. On the other hand, $\tilde{\mathcal{B}}$ is injective (i.e., $\tilde{\mathcal{A}}_{\a}(S_\ga)$ does not contain nontrivial flat functions, and then the class $\tilde{\mathcal{A}}_{\a}(S_\ga)$ is said to be quasianalytic) if, and only if, $\ga>\a$ (Watson's lemma, see for example~\cite[Prop.\ 11]{balserutx}). Our main aim in this paper is to provide generalizations of this kind of results in the framework of Carleman ultraholomorphic classes associated with strongly regular sequences inducing a proximate order. Let us start with an overview of the existing literature in this respect.\par

In 1995 V. Thilliez~\cite{Thilliez1} obtained right inverses in the Gevrey case when $\ga<\a$ by applying techniques from the ultradifferentiable setting (i.e. regarding extension results for classes of smooth functions on open subsets of $\mathbb{R}^n$, determined by imposing a suitable growth of the derivatives), and the same was done by the author in~\cite{javier} by adapting the truncated Laplace transform procedure already used by J. P. Ramis in Borel--Ritt--Gevrey theorem~\cite{Ramis1} (this second solution was well-suited for the extension of this result to the several variable case). Regarding general classes, J. Schmets and M. Valdivia~\cite{SchmetsValdivia} extended some results of H.-J. Petzsche~\cite{pet} for ultradifferentiable classes, and applied them in order to provide the first powerful results in the present framework. Subsequently, V. Thilliez~\cite{thilliez} improved the results in~\cite{SchmetsValdivia} in several respects (see Subsection 3.1 in his paper for the details) by relying on a double application of suitable Whitney's extension results for Whitney ultradifferentiable jets on compact sets with Lipschitz boundary appearing in \cite{pet,bbmt,chaucho}. In particular, he introduced a growth index $\ga(\M)\in(0,\infty)$ for every strongly regular sequence $\M$ (which for $\M_{\a}$ equals $\a$), and proved the following facts: if $\ga<\ga(\M)$, then $\mathcal{A}_{\M}(S_{\ga})$ is not quasianalytic, and there exist right inverses for $\tilde{\mathcal{B}}$, which are obtained due to the explicit construction of nontrivial flat functions in the class $\mathcal{A}_{\M}(S_\ga)$. Indeed, these flat functions allowed A. Lastra, S. Malek and the author~\cite{lastramaleksanz} to define suitable kernels and moment sequences by means of which to obtain again right inverses by the classical truncated Laplace transform technique. Because of the integral form of the solution, this procedure admits an easy generalization to the several variable case, and does not rest on any result from the ultradifferentiable setting.\par
However, the preceding results for general classes are not fully satisfactory. Firstly, the equivalences stated in Borel--Ritt--Gevrey theorem and Watson's lemma for the Gevrey case are now only one-side implications. Secondly, and strongly related to the previous remark, the need to restrict the opening of the sector $S_\ga$ to $\ga<\ga(\M)$ in order to obtain flat functions in $\mathcal{A}_{\M}(S_{\ga})$ does not allow one to treat the apparently limit situation in which $\ga=\ga(\M)$. Note that, in the Gevrey case, the function $e^{-z^{-1/\a}}$ is flat in the class $\tilde{\mathcal{A}}_{\a}(S_{\a})$, and of course in every  $\mathcal{A}_{\a}(S_{\ga})$ for $\ga<\a$. So, our main objective will be to obtain flat functions in sectors of optimal opening.\par

In this sense, we first introduce for every strongly regular sequence $\M$ a new constant $\omega(\M)$, measuring the rate of growth of the sequence $\M$, in terms of which quasianalyticity in the classes $\mathcal{A}_{\M}(S_{\ga})$ may be properly characterized due to a classical result of B. I. Korenbljum (\cite{korenbljum}; see Theorem~\ref{teorKoren}). This constant is easily computed in concrete situations (see~(\ref{equaordequasM})), and indeed it is the inverse of the order of growth of the classical function $M(t)$ associated with $\M$, namely
$M(t)=\sup_{n\in\N_0}\log(t^n/M_n)$, $t>0$ (see~(\ref{equaordequasMbis})). Regarding the construction of flat functions, V. Thilliez (\cite{Thilliez2}) had characterized flatness in $\tilde{\mathcal{A}}_{\M}(S_{\ga})$ in terms of the existence of non-uniform estimates governed by the function $e^{-M(1/|z|)}$, much in the same way as the function $e^{-z^{-1/\a}}$ expresses flatness in the Gevrey case.
So, it became clear to us the need to construct functions in sectors whose growth is accurately specified by the function $M(t)$. The classical theory of growth for holomorphic functions defined in sectorial regions, based on the notion of (constant) exponential order, showed itself not profound enough to deal with the general case. Luckily, the theory of proximate orders, allowing to change the constant order $\rho>0$ into a function $\rho(r)$ more closely specifying the desired rate of growth, is available since the 1920s, and some of its quite recent developments, mainly due to L. S. Maergoiz~\cite{Maergoiz}, have been the key for our success. The problem of characterizing those sequences $\M$ associated with a proximate order has been solved (Proposition~\ref{propcaracdderordenaprox}), and it turns out that all the interesting examples of strongly regular sequences appearing in the literature belong to this class.
Whenever this is the case, the results of L. S. Maergoiz allow us to obtain the desired flat functions in $\tilde{\mathcal{A}}_{\M}(S_{\omega(\M)})$ (see Theorem~\ref{teorconstrfuncplana}) and, immediately, we may generalize Watson's lemma, see Corollary~\ref{coroWatsonlemma}. Subsequently, in Section~\ref{sectkernels} suitable kernels and moment sequences are introduced, by means of which we may prove that $\tilde{\mathcal{B}}$ is surjective in $\tilde{\mathcal{A}}_{\M}(S_{\ga})$ if, and only if, $\ga\le\omega(\M)$, so generalizing Borel--Ritt--Gevrey theorem (see Theorem~\ref{tpral}).\par

It should be mentioned that for the standard strongly regular sequences appearing in the literature, the value of the constants $\ga(\M)$ and $\omega(\M)$ agree. However, we have only been able to prove that $\ga(\M)\le \omega(\M)$ in general. In case $\ga(\M)<\omega(\M)$ can actually occur for some sequences, our results would definitely improve those of V. Thilliez by enlarging the sectors for which non-quasianalyticity holds or right inverses exist. In any case, the equivalences stated in Theorem \ref{tpral} and Corollary \ref{coroWatsonlemma} are new.\par

Finally, in Theorem~\ref{teorinversasderecha} we gather the information concerning the existence of right inverses for $\tilde{\mathcal{B}}$ in $\mathcal{A}_{\M}(S_{\ga})$: they exist whenever $\ga<\omega(\M)$, and their existence, under some specific condition (satisfied, for instance, in the Gevrey case), implies that $\ga<\omega(\M)$.

\section{Preliminaries}
\subsection{Notation}
We set $\N:=\{1,2,...\}$, $\N_{0}:=\N\cup\{0\}$.
$\mathcal{R}$ stands for the Riemann surface of the logarithm, and
$\C[[z]]$ is the space of formal power series in $z$ with complex coefficients.\par\noindent
For $\gamma>0$, we consider unbounded sectors
$$S_{\gamma}:=\{z\in\mathcal{R}:|\hbox{arg}(z)|<\frac{\gamma\,\pi}{2}\}$$
or, in general, bounded or unbounded sectors
$$S(d,\alpha,r):=\{z\in\mathcal{R}:|\hbox{arg}(z)-d|<\frac{\alpha\,\pi}{2},\ |z|<r\},\quad
S(d,\alpha):=\{z\in\mathcal{R}:|\hbox{arg}(z)-d|<\frac{\alpha\,\pi}{2}\}$$
with bisecting direction $d\in\R$, opening $\alpha\,\pi$ and (in the first case) radius $r\in(0,\infty)$.\par\noindent
A sectorial region $G(d,\a)$ with bisecting direction $d\in\R$ and opening $\alpha\,\pi$ will be a domain in $\mathcal{R}$ such that $G(d,\a)\subset S(d,\a)$, and
for every $\beta\in(0,\a)$ there exists $\rho=\rho(\beta)>0$ with $S(d,\beta,\rho)\subset G(d,\a)$. In particular, sectors are sectorial regions.\par\noindent
A sector $T$ is a bounded proper subsector of a sectorial region $G$ (denoted by $T\ll G$) whenever the radius of $T$ is finite and $\overline{T}\setminus\{0\}\subset G$.
Given two unbounded sectors $T$ and $S$, we say $T$ is an unbounded proper subsector of $S$, and we write $T\prec S$, if $\overline{T}\setminus\{0\}\subset S$.\par\noindent
$\mathcal{H}(U)$ denotes the space of holomorphic functions in an open set $U\subset\mathcal{R}$.\par\noindent
$D(z_0,r)$ stands for the disk centered at $z_0$ with radius $r>0$.

\subsection{Asymptotic expansions and ultraholomorphic classes}\label{subsectCarlemanclasses}

Given a sequence of positive real numbers $\M=(M_n)_{n\in\N_0}$, a constant $A>0$ and a sector $S$, we define
$$\mathcal{A}_{\M,A}(S)=\big\{f\in\mathcal{H}(S):\left\|f\right\|_{\M,A}:=\sup_{z\in S,n\in\N_{0}}\frac{|f^{(n)}(z)|}{A^{n}n!M_{n}}<\infty\big\}.$$
%\begin{itemize}
%\item $f\in\mathcal{A}_{\M,A}(S): |f^{(n)}(z)|\le\left\|f\right\|_{\M,A}A^{n}n!M_n$, $n\in\N_{0},z\in S$.
%\item
($\mathcal{A}_{\M,A}(S),\left\| \ \right\| _{\M,A}$) is a Banach space, and $\mathcal{A}_{\M}(S):=\cup_{A>0}\mathcal{A}_{\M,A}(S)$ is called a \textit{Carleman ultraholomorphic class} in the sector $S$.
\parn
One may accordingly define classes of sequences
$$\Lambda_{\M,A}=\Big\{\bmu=(\mu_{n})_{n\in\N_{0}}\in\C^{\N_{0}}: \left|\bmu\right|_{\M,A}:=\sup_{n\in\N_{0}}\displaystyle \frac{|\mu_{n}|}{A^{n}n!M_{n}}<\infty\Big\}.$$
$(\Lambda_{\M,A},\left| \  \right|_{\M,A})$ is again a Banach space, and we put $\Lambda_{\M}:=\cup_{A>0}\Lambda_{\M,A}$.\parn
%\end{itemize}
Since the derivatives of $f\in\mathcal{A}_{\bM,A}(S)$ are Lipschitzian, for every $n\in\N_{0}$ one may define
$$f^{(n)}(0):=\lim_{z\in S,z\to0 }f^{(n)}(z)\in\C,$$
and it is clear that the sequence
\begin{equation*}%\label{equaBoremap}
\tilde{\mathcal{B}}(f):=(f^{(n)}(0))_{n\in\N_{0}}\in\Lambda_{\M,A},\qquad f\in\mathcal{A}_{\bM,A}(S).
\end{equation*}
The map $\tilde{\mathcal{B}}:\mathcal{A}_{\M}(S)\longrightarrow \Lambda_{\M}$
so defined is the \textit{asymptotic Borel map}.

Next, we will recall the relationship between these classes and the concept of asymptotic expansion.
\begin{defi}
We say a holomorphic function $f$ in a sectorial region $G$ admits the formal power series $\hat{f}=\sum_{p=0}^{\infty}a_{p}z^{p}\in\C[[z]]$ as its $\bM-$\emph{asymptotic expansion} in $G$ (when the variable tends to 0) if for every $T\ll G$ there exist $C_T,A_T>0$ such that for every $n\in\N$, one has
\begin{equation*}\Big|f(z)-\sum_{p=0}^{n-1}a_pz^p \Big|\le C_TA_T^nM_{n}|z|^n,\qquad z\in T.%\label{desarasint}
\end{equation*}
We will write $f\sim_{\bM}\sum_{p=0}^{\infty}a_pz^p$ in $G$. $\tilde{\mathcal{A}}_{\M}(G)$ stands for the space of functions admitting $\bM-$asymptotic expansion in $G$.\par
\end{defi}

\begin{defi}
Given a sector $S$, we say $f\in\mathcal{H}(S)$ admits $\hat{f}$ as its \emph{uniform $\bM-$asymptotic expansion in $S$ of type $A>0$} if there exists $C>0$ such that for every $n\in\N$, one has
\begin{equation}\Big|f(z)-\sum_{p=0}^{n-1}a_pz^p \Big|\le CA^nM_{n}|z|^n,\qquad z\in S.\label{desarasintunifo}
\end{equation}
\end{defi}

As a consequence of Taylor's formula and Cauchy's integral formula for the derivatives, we have the following result (see \cite{balserutx,galindosanz}).
\begin{prop}\label{propcotaderidesaasin}
Let $S$ be a sector and $G$ a sectorial region.
\begin{itemize}
\item[(i)] If $f\in\mathcal{A}_{\M,A}(S)$, then $f$ admits $\hat{f}=\sum_{p\in\N_0}\frac{1}{p!}f^{(p)}(0)z^p$ as its uniform $\bM-$asymptotic expansion in $S$ of type $A$.
\item[(ii)] $f\in\tilde{\mathcal{A}}_{\M}(G)$ if, and only if, for every $T\ll G$ there exists $A_T>0$ such that $f|_T\in \mathcal{A}_{\M,A_T}(T)$. Hence, the map $\tilde{\mathcal{B}}:\tilde{\mathcal{A}}_{\M}(G)\longrightarrow \Lambda_{\M}$ is also well defined.
\end{itemize}
\end{prop}

\begin{defi}
A function $f$ in any of the previous classes is said to be \textit{flat} if $\tilde{\mathcal{B}}(f)$ is the null sequence or, in other words, $f\sim_{\M}\hat{0}$, where $\hat{0}$ denotes the null power series.
\end{defi}

\begin{rema}\label{remaCarlclassasympexpan}
As a consequence of Cauchy's integral formula for the derivatives, given a sector $S$ one can prove that whenever $T\ll S$, there exists a constant $c=c(T,S)>0$ such that the restriction to $T$, $f_T$, of functions $f$ defined on $S$ and admitting uniform $\bM-$asymptotic expansion in $S$ of type $A>0$, belongs to $\mathcal{A}_{\bM,cA}(T)$, and moreover, if one has (\ref{desarasintunifo}) then $\Vert f_T\Vert_{\bM,cA}\le C$.
\end{rema}

\subsection{Strongly regular sequences and associated functions}\label{subsectstrregseq}

Most of the information in this subsection is taken from the works of A. A. Goldberg and I. V. Ostrovskii~\cite{GoldbergOstrowskii}, H. Komatsu~\cite{komatsu} and V. Thilliez~\cite{thilliez}, which we refer to for further details and proofs.
In what follows, $\bM=(M_p)_{p\in\N_0}$ will always stand for a sequence of
positive real numbers, and we will always assume that $M_0=1$.
\begin{defi}\label{sucfuereg}
We say $\bM$ is \textit{strongly regular} if the following hold:\par
($\a_0$) $\bM$ is \textit{logarithmically convex}: $M_{p}^{2}\le M_{p-1}M_{p+1}$ for
every $p\in\N$.\par
($\mu$) $\bM$ is of \textit{moderate growth}: there exists $A>0$ such that
$$M_{p+\ell}\le A^{p+\ell}M_{p}M_{\ell},\qquad p,\ell\in\N_0.$$
\par($\gamma_1$) $\bM$ satisfies the \textit{strong non-quasianalyticity condition}: there exists $B>0$ such that
$$
\sum_{\ell\ge p}\frac{M_{\ell}}{(\ell+1)M_{\ell+1}}\le B\frac{M_{p}}{M_{p+1}},\qquad p\in\N_0.
$$
\end{defi}

\begin{rema}
In the literature a different set of conditions appears frequently when dealing with ultraholomorphic or ultradifferentiable classes of functions. Let us clarify the relationship between these two approaches:
If $\M$ is strongly regular, then $\M'=(n!M_n)_{n\in\N_0}$ verifies the standard conditions (M.1), (M.2) and (M.3) of H. Komatsu (see \cite{komatsu,mandelbrojt}). On the other hand,
if a sequence of positive real numbers $\M'=(M'_n)_{n\in\N_0}$, with $M'_0=1$, verifies (M.2) and (M.3) of H. Komatsu, and moreover $\M:=(M'_n/n!)_{n\in\N_0}$ is logarithmically convex, then $\M$ is strongly regular.
\end{rema}

\begin{exam}\label{examstroregusequ}
\begin{itemize}
\item[(i)] The best known example of strongly regular sequence is $\M_{\a}=(n!^{\a})_{n\in\N_{0}}$, called the \textit{Gevrey sequence of order $\a>0$}.
\item[(ii)] The sequences $\M_{\a,\b}=\big(n!^{\a}\prod_{m=0}^n\log^{\b}(e+m)\big)_{n\in\N_0}$, where $\a>0$ and $\b\in\R$, are strongly regular.
\item[(iii)] For $q>1$, $\M=(q^{n^2})_{n\in\N_0}$ is logarithmically convex and strongly non-quasianalytic, but not of moderate growth.
\end{itemize}
\end{exam}

For a sequence $\bM=(M_{p})_{p\in\N_0}$ verifying properties $(\alpha_0)$ and $(\gamma_{1})$ one has that the associated \textit{sequence of quotients}, $\bm=(m_{p}:=M_{p+1}/M_{p})_{p\in\N_0}$, is an increasing sequence to infinity, so that the map $h_{\bM}:[0,\infty)\to\R$, defined by
\begin{equation*}%\label{equadefihdeM}
h_{\bM}(t):=\inf_{p\in\N_{0}}M_{p}t^p,\quad t>0;\qquad h_{\bM}(0)=0,
\end{equation*}
turns out to be a non-decreasing continuous map in $[0,\infty)$ onto $[0,1]$. In fact
$$
h_{\bM}(t)= \left \{ \begin{matrix}  t^{p}M_{p} & \mbox{if }t\in\big[\frac{1}{m_{p}},\frac{1}{m_{p-1}}\big),\ p=1,2,\ldots,\\
1 & \mbox{if } t\ge 1/m_{0}. \end{matrix}\right.
$$

\begin{defi}[\cite{pet},~\cite{chaucho}]
Two sequences $\bM=(M_{p})_{p\in\N_0}$ and $\bM'=(M'_{p})_{p\in\N_0}$ of positive real numbers are said to be \textit{equivalent} if there exist positive constants $L,H$ such that
$$L^pM_p\le M'_p\le H^pM_p,\qquad p\in\N_0.$$
\end{defi}
In this case, it is straightforward to check that
\begin{equation}\label{equahdeMequi}
h_{\bM}(Lt)\le h_{\bM'}(t)\le h_{\bM}(Ht),\qquad t\ge 0.
\end{equation}

One may also associate with a strongly regular sequence $\bM$ the function
\begin{equation}\label{equadefiMdet}
M(t):=\sup_{p\in\N_{0}}\log\big(\frac{t^p}{M_{p}}\big)=-\log\big(h_{\bM}(1/t)\big),\quad t>0;\qquad M(0)=0,
\end{equation}
which is a non-decreasing continuous map in $[0,\infty)$ with $\lim_{t\to\infty}M(t)=\infty$.
Indeed,
\begin{equation*}%\label{equaexprMdet}
M(t)=\left \{ \begin{matrix}  p\log t -\log(M_{p}) & \mbox{if }t\in [m_{p-1},m_{p}),\ p=1,2,\ldots,\\
0 & \mbox{if } t\in [0,m_{0}), \end{matrix}\right.
\end{equation*}
and one can easily check that $M$ is convex in $\log t$, i.e., the map $t\mapsto M(e^t)$ is convex in $\mathbb{R}$.\parn

Some additional properties of strongly regular sequences needed in the present work are the following ones.
\begin{lemma}[\cite{thilliez}]
Let $\M=(M_{p})_{p\in\N_0}$ be a strongly regular sequence and $A>0$ the constant appearing in the property $(\mu)$ in Definition~\ref{sucfuereg}. Then,
\begin{equation}
\label{e107}
m_{p}\le A^{2}M_{p}^{1/p}\le A^{2}m_{p}\qquad \hbox{for every }p\in\N_0.
\end{equation}
Let $s$ be a real number with $s\ge1$. There exists $\rho(s)\ge1$ (only depending on $s$ and $\M$) such that
\begin{equation}\label{e120}
h_{\M}(t)\le(h_{\M}(\rho(s)t))^{s}\qquad\hbox{for }t\ge0.
\end{equation}
\end{lemma}

\begin{rema}\label{remapropsuceemepealape}
\begin{itemize}
\item[(i)] The condition of moderate growth $(\mu)$ plays a fundamental role in the proof of (\ref{e120}), which will in turn be crucial in many of our arguments.
\item[(ii)] From property (\ref{e107}) we deduce that $\bM$ and $(m_{p}^p)_{p\in\N_0}$ are equivalent.
\item[(iii)] For every $p\in\N_0$, the continuity of $M$ at $m_p$ amounts to the trivial equality $m_{p}^p/M_p=m_{p}^{p+1}/M_{p+1}$.
\item[(iv)] Moreover, since the sequence $\bm=(m_p)_{p\in\N_0}$ (respectively, the function $M(t)$) increases to infinity as $p$ (resp. $t$) tends to infinity, the sequence $(M(m_p))_{p\in\N_0}=\big(\log(m_{p}^p/M_p)\big)_{p\in\N_0}$, and consequently also $(m_{p}^p/M_p)_{p\in\N_0}$, increase to infinity, starting at the value 0 and 1, respectively. Note also that the $p$-th and $(p+1)$-th terms of any of these two sequences are equal if, and only if, $m_p=m_{p+1}$.
\end{itemize}
\end{rema}

We now recall the following definitions and facts, mainly taken from the book of A. A. Goldberg and I. V. Ostrovskii~\cite{GoldbergOstrowskii}.

\begin{defi}[\cite{GoldbergOstrowskii}, p.\ 43]%\label{defiordefunc}
Let $\a(r)$ be a nonnegative and nondecreasing function in $(c,\infty)$ for some $c\ge 0$ (we write $\a\in\Lambda$). The \textit{order} of $\a$ is defined as
$$
\rho=\rho[\a]:=\limsup_{r\to\infty}\frac{\log^+\a(r)}{\log r}\in[0,\infty],
$$
where $\log^+=\max(\log,0)$. $\a(r)$ is said to have finite order if $\rho<\infty$.
\end{defi}

We are firstly interested in determining the order of the function $M(r)\in\Lambda$ defined in (\ref{equadefiMdet}) and associated with a strongly regular sequence $\M$.
To this end, we need to recall now the definition of exponent of convergence of a sequence and how it may be computed.

\begin{prop}[\cite{holland}, p.\ 65]
Let $(c_n)_{n\in\N_0}$ be a nondecreasing sequence of positive real numbers tending to infinity. The \textit{exponent of convergence} of $(c_n)_n$ is defined as
$$
\lambda_{(c_n)}:=\inf\{\mu>0:\sum_{n=0}^\infty \frac{1}{c_n^{\mu}}\textrm{ converges}\}
$$
(if the previous set is empty, we put $\lambda_{(c_n)}=\infty$). Then, one has
\begin{equation}\label{equaexpoconv}
\lambda_{(c_n)}=\limsup_{n\to\infty}\frac{\log(n)}{\log(c_n)}.
\end{equation}
\end{prop}

We will also need the following fact, which can be found in~\cite{mandelbrojt}: if we consider the \textit{counting function} for the sequence of quotients $\bm$,
$\nu:(0,\infty)\to\N_0$ given by
\begin{equation}\label{equadefinuder}
\nu(r):=\#\{j:m_j\le r\},
\end{equation}
then one has that
\begin{equation}\label{equarelaMdetnuder}
M(t)=\int_0^t\frac{\nu(r)}{r}\,dr,\qquad t>0.
\end{equation}

We may now state our first result.
\begin{theo}\label{teorordeMdetiguallimidder}
Let $\M$ be strongly regular, $\bm$ the sequence of its quotients and $M(r)$ its associated function. Then, the order of $M(r)$ is given by
\begin{equation}\label{equaordeMdet}
\rho[M]=\lim_{r\to\infty}\frac{\log M(r)}{\log r}=\limsup_{n\to\infty}\frac{\log(n)}{\log(m_{n})}.
\end{equation}
\end{theo}

\begin{rema}\label{remadefifuncdder}
The function $d$, defined for $r>\max\{1,m_0\}$ by $d(r)=\log(M(r))/\log r$, will play an important role in what follows. It is clearly continuous and piecewise continuously differentiable in its domain (meaning that it is differentiable except at a sequence of points, tending to infinity, at any of which it is continuous and has distinct finite lateral derivatives).
\end{rema}

\begin{proof1}
The first equality for $\rho[M]$ is due to the fact that the function $d(r)$ is eventually strictly increasing, as we now show. It is enough to prove that $d'(r)> 0$ for $r\in(m_{p-1},m_p)$ and $p$ large enough. This is best seen by considering the auxiliary function
$$
D(t):=d(e^t)=\frac{\log\big(pt-\log(M_p)\big)}{t},\qquad t\in(\log(m_{p-1}),\log(m_p)),\ p\in\N,
$$
and then proving that $D'(t)> 0$ for $t\in(\log(m_{p-1}),\log(m_p))$ and large enough. We have
$$
D'(t)=\frac{1}{t^2}\Big(1+\frac{\log(M_p)}{pt-\log(M_p)}-\log\big(pt-\log(M_p)\big)\Big),\qquad t\in(\log(m_{p-1}),\log(m_p)).
$$
When $t$ runs over $(\log(m_{p-1}),\log(m_p))$, the value $pt-\log(M_p)$ runs over
$$
\big(\log(m_{p-1}^{p-1}/M_{p-1}),\log(m_{p}^p/M_p)\big),
$$
which, as long as $m_{p-1}<m_p$, is a nonempty interval contained in $(0,\infty)$ (see Remark~\ref{remapropsuceemepealape}.(iv)).\par
Observe that $M_p>1$ for $p$ large enough, what we assume from now on. Since for any $A>0$ the function $y\in(0,\infty)\mapsto 1+A/y-\log(y)$ is strictly decreasing, we will conclude that $D'(t)> 0$ for $t\in(\log(m_{p-1}),\log(m_p))$ if we have that
$$
1+\frac{\log(M_p)}{\log(m_p^p/M_p)}-\log(\log(m_p^p/M_p))> 0.
$$
But
$$
1+\frac{\log(M_p)}{\log(m_p^p/M_p)}-\log(\log(m_p^p/M_p))=
\frac{\log\Big(\frac{m_p^p}{\log(m_p^p/M_p)}\Big)}{\log(m_p^p/M_p)},
$$
whose denominator is positive; finally, note that
$$
m_p^p>\log(m_p^p)>\log(m_p^p)-\log(M_p)=\log(m_p^p/M_p)>0,
$$
so that the numerator is also positive and we are done.\par
For the second expression of $\rho[M]$, we take into account the link given in (\ref{equarelaMdetnuder}) between $M(r)$ and  the counting function $\nu(r)$ for $\bm$ (as defined in (\ref{equadefinuder})), which also belongs to $\Lambda$. We may apply Theorem 2.1.1 in~\cite{GoldbergOstrowskii} to deduce that the order of $M(r)$ equals that of $\nu(r)$.
Now, from Theorem 2.1.8 in~\cite{GoldbergOstrowskii} we know that the order of $\nu(r)$ is in turn the exponent of convergence of $\bm$, given by the formula in (\ref{equaexpoconv}).
\end{proof1}

\begin{rema}\label{remaproprhoM}
\begin{itemize}
\item[(i)] Let $\M_{\a}$ be the Gevrey sequence of order $\a>0$, and $M_{\a}(r)$ its associated function. By means of (\ref{equaordeMdet}), it is obvious that $\rho[M_{\a}]=1/\a$.
The same is true for any sequence of the form $(a^nn!^{\a})_{n\in\N_0}$, with $a,\a>0$.
\item[(ii)] Let $\M=(M_n)_{n\in\N_0}$ and $\M^*=(M^{*}_n)_{n\in\N_0}$ be strongly regular sequences such that $M_n\le M^{*}_n$ for every $n\in\N_0$. By the very definition of the respective associated functions $M(r)$ and $M^{*}(r)$, one has that $M(r)\ge M^{*}(r)$ for every $r\ge 0$, and consequently $\rho[M]\ge \rho[M^{*}]$.
\item[(iii)] By Lemma 1.3.2 in Thilliez~\cite{thilliez}, for every strongly regular sequence $\M$ there exist positive constants $a_1,a_2,\ga,\delta$, with $\ga<\delta$, such that
$$
a_1^nn!^{\ga}\le M_n\le a_2^nn!^{\delta},\qquad n\in\N_0.
$$
From the two previous remarks we deduce that $1/\delta\le\rho[M]\le 1/\ga$, and, in particular, $\rho[M]\in(0,\infty)$.
\end{itemize}
\end{rema}

We next recall the notion of growth index defined and studied by V. Thilliez~\cite[Sect.\ 1.3]{thilliez}.

\begin{defi}\label{defiindegrowM}
Let $\bM=(M_{p})_{p\in\N_{0}}$ be a strongly regular sequence and $\ga>0$. We say $\bM$ satisfies property $\left(P_{\ga}\right)$  if there exist a sequence of real numbers $m'=(m'_{p})_{p\in\N_0}$ and a constant $a\ge1$ such that: (i) $a^{-1}m_{p}\le m'_{p}\le am_{p}$, $p\in\N$, and (ii) $\left((p+1)^{-\ga}m'_{p}\right)_{p\in\N_0}$ is increasing.

The \textit{growth index} of $\bM$ is
$$\ga(\bM):=\sup\{\ga\in\R:(P_{\ga})\hbox{ is fulfilled}\}\in(0,\infty).$$
\end{defi}

\begin{exam}\label{examvalogamm}
\begin{itemize}
\item[(i)] For the Gevrey sequence of order $\a>0$, one has $\ga(\M_{\a})=\a$.
\item[(ii)] For the sequences $\M_{\a,\b}$ in Example~\ref{examstroregusequ}.(ii) one can check that $\gamma(\M_{\a,\b})=\a$.
%\item[(iii)] For $q>1$, $\M=(q^{n^2})_{n\in\N_0}$ is logarithmically convex and strongly non-quasianalytic, but not of moderate growth.
\end{itemize}
\end{exam}

\section{Results on quasianalyticity in ultraholomorphic classes}%\label{sectquasianal}

We are interested in characterizing those classes in which the asymptotic Borel map is injective. First, quasianalytic Carleman classes are defined.

\begin{defi}
Let $S$ be a sector and $\M=(M_{p})_{p\in\N_{0}}$ be a sequence of positive numbers. We say that $\mathcal{A}_{\M}(S)$ is  \textit{quasianalytic} if it does not contain nontrivial flat functions.
\end{defi}\parn

Characterizations of quasianalyticity for general sequences $\M$ in one and several variables are available in \cite{lastrasanz1}, generalizing the work of B. I. Korenbljum \cite{korenbljum}. In this paper, we restrict our attention to the one-variable case, and focus on strongly regular sequences, although in many of the results in this section weaker assumptions on $\M$ suffice. As shown in the next result, quasianalyticity is governed by the opening of the sector.

\begin{theo}[\cite{korenbljum}]\label{teorKoren}
Let $\M$ be strongly regular and $\ga>0$. The following statements are equivalent:
\begin{itemize}
\item[(i)] The class $\mathcal{A}_{\M}(S_{\gamma})$ is quasianalytic.
\item[(ii)] $\displaystyle \sum_{n=0}^{\infty}\Big(\frac{M_{n}}{(n+1)M_{n+1}}\Big)^{1/(\ga+1)}=\infty$.
\end{itemize}
\end{theo}

Accordingly, we introduce a new quantity regarding quasianalyticity.
\begin{defi}%\label{defiordequas}
For a strongly regular sequence $\M$, we put
$$Q_{\M}=\{\ga>0:\mathcal{A}_{\M}(S_{\gamma}) \textrm{ is quasianalytic}\}.
$$
The \textit{order of quasianalyticity} of $\M$ is defined as $\omega(\M):=\inf Q_{\M}$.
\end{defi}

We can obtain its value due to the following result.
\begin{theo}\label{teororderM}
For a strongly regular sequence $\M$ with associated function $M(r)$, we have
\begin{equation}\label{equaordequasM}
\omega(\M)=\liminf_{n\to\infty}
\frac{\log(m_{n})}{\log(n)}=\frac{1}{\lambda_{(m_n)}},
\end{equation}
and consequently,
\begin{equation}\label{equaordequasMbis}
\omega(\M)=\frac{1}{\rho[M]}\in(0,\infty).
\end{equation}
\end{theo}

\begin{proof1}
Since $\M$ is strongly regular, $(n!M_n)_{n\in\N_0}$ is logarithmically convex. So, the sequence of its quotients, $(\frac{(n+1)M_{n+1}}{M_{n}})_{n\in\N_0}$, is nondecreasing and, moreover,
tends to infinity because of property $(\gamma_1)$ in Definition~\ref{sucfuereg}. In view of (\ref{equaexpoconv}), the exponent of convergence of the sequence $((n+1)M_{n+1}/M_{n})_{n\in\N_0}=((n+1)m_{n})_{n\in\N_0}$ and that of the sequence $(m_{n})_{n\in\N_0}$ are related as follows:
\begin{equation*}
\lambda_{((n+1)m_{n})}= \limsup_{n\to\infty}\frac{\log(n)}{\log((n+1)m_{n})}
=\frac{1}{1+\liminf_{n\to\infty}\frac{\log(m_n)}{\log(n)}}=\frac{1}{1+1/\lambda_{(m_n)}}.
\end{equation*}
On the other hand, from Theorem~\ref{teorKoren} and the definition of $\omega(\M)$ it is clear that $$\frac{1}{\omega(\M)+1}=\lambda_{((n+1)m_{n})},
$$
hence
\begin{equation*}%\label{equaomegM}
\omega(\M)=\liminf_{n\to\infty}\frac{\log(m_{n})}{\log(n)}=\frac{1}{\lambda_{(m_n)}}.
\end{equation*}
Comparing this to (\ref{equaordeMdet}), and by Remark~\ref{remaproprhoM}.(iii), we conclude.
\end{proof1}

\begin{rema}%\label{remaoptiopen}
Observe that $\pi\omega(\M)$ is the optimal opening for quasianalyticity, in the sense that
the class $\mathcal{A}_{\M}(S)$ is (respectively, is not) quasianalytic whenever the opening of $S$ exceeds (resp. is less than) this quantity. When the opening of the sector equals $\pi\omega(\M)$ both cases are possible, as shown in the forthcoming Example~\ref{examsucesalfabetaomegagamma}.
\end{rema}

\begin{rema}%\label{remapropomegM}
\begin{itemize}
\item[(i)] Consider a pair of equivalent sequences $\bM$ and $\bM'$. Given a sector $S$, the spaces $\mathcal{A}_{\bM}(S)$ and $\mathcal{A}_{\bM'}(S)$ coincide. For a sectorial region $G$, also $\tilde{\mathcal{A}}_{\M}(G)$ and $\tilde{\mathcal{A}}_{\M'}(G)$ agree. So, it is clear that $\omega(\M)=\omega(\M')$, and from (\ref{equaordequasMbis}), also $\rho[M]=\rho[M']$ for the associated functions (this last fact can also be deduced from~(\ref{equahdeMequi})).
\item[(ii)] If the strongly regular sequences $\M$ and $\M^{*}$ are such that $M_n\le M^{*}_n$ for every $n$, then $\mathcal{A}_{\bM}(S)\subset\mathcal{A}_{\bM^{*}}(S)$ for any sector $S$, and so $Q_{\M^{*}}\subset Q_{\M}$ and $\omega(\M)\le \omega(\M^{*})$. Note that this fact is not at all clear from the formula~(\ref{equaordequasM}).
\end{itemize}
\end{rema}

Regarding the relationship between $\omega(\M)$ and $\ga(\M)$, we have the following result.
\begin{prop}%\label{proprelagammomegM}
For every strongly regular sequence $\M$ one has $\omega(\M)\ge \gamma(\M)$.
\end{prop}
\begin{proof1}
Suppose $\M$ verifies $(P_\ga)$ (see Definition~\ref{defiindegrowM}) for some $\ga>0$.
As indicated in~\cite[p.~173]{thilliez}, this easily
implies the existence of a constant $a>0$ such that
$a^nn!^{\ga}\le M_n$ for every $n\in\N_0$. Hence, by (ii) in the previous remark we have
$\ga=\omega((a^nn!^{\ga})_{n\in\N_0})\le\omega(\M)$, and the definition of $\ga(\M)$ is enough to conclude.
\end{proof1}

According to the very definition of $\omega(\M)$, the previous result is indeed equivalent to the following one, proved by V. Thilliez \cite{thilliez} and, subsequently, by A. Lastra and the author \cite{lastrasanz1}.
However, the present argument seems to be simpler than the ones involved in the previous proofs of this theorem.

\begin{theo}[\cite{thilliez,lastrasanz1}]
Let $0<\ga<\ga(\M)$.
Then, the class $\mathcal{A}_{\M}(S_{\gamma})$ is not quasianalytic.
\end{theo}

It is an open problem to decide whether $\omega(\M)=\ga(\M)$ in general. At the moment, we have not been able to find an example showing that $\omega(\M)>\gamma(\M)$ may occur.
However, from the fact that $\omega(\M)=1/\lambda_{(m_n)}$ we may deduce an easy characterization for the equality of both constants.

\begin{coro}%[\cite{lastrasanz1}, Thm.\ 4.10]\label{teorWatson}
Let $\M$ be strongly regular. The following statements are equivalent:
\begin{itemize}
\item[(i)] $\omega(\M)=\ga(\M)$,
\item[(ii)] For every $\gamma>\gamma(\M)$ we have that
\begin{equation*}
\sum_{n=0}^{\infty}\big(\frac{1}{m_{n}}\big)^{1/\ga} =\infty.%\label{conditionWatson}
\end{equation*}
\end{itemize}
\end{coro}

\begin{exam}\label{examsucesalfabetaomegagamma}
Consider the sequences $\M_{\a,\b}$, $\a>0$, $\b\in\R$, introduced in Examples~\ref{examstroregusequ}.(ii) and~\ref{examvalogamm}.(ii). Applying Theorem \ref{teorKoren}, it is easy to check that
$$Q_{\M_{\a,\b}}=\begin{cases}[\a,\infty)&\textrm{ if }\a\ge \b-1,\\
(\a,\infty)&\textrm{ if }\a< \b-1,\end{cases}
$$
so that $\omega(\M_{\a,\b})=\a=\gamma(\M_{\a,\b})$.
\end{exam}

\section{Flat functions via proximate orders}%\label{sectflatfunctions}

In this section we show how one can construct flat functions in the classes $\tilde{\mathcal{A}}_{\M}(S_{\omega(\M)})$ for a strongly regular sequence $\M$ by relying on the notion of analytic proximate orders, appearing in the theory of growth of entire functions and developed, among others, by E. Lindel\"of, G. Valiron, B. Ja. Levin, A. A. Goldberg, I. V. Ostrosvkii and L. S. Maergoiz (see the references \cite{Valiron,Levin,GoldbergOstrowskii,Maergoiz}).

\begin{defi}[\cite{Valiron}]%\label{defiproxorde}
We say a real function $\ro(r)$, defined on $(c,\infty)$ for some $c\ge 0$, is a \textit{proximate order} if the following hold:
 \begin{enumerate}[(i)]
  \item $\rho(r)$ is continuous and piecewise continuously differentiable in $(c,\infty)$,
  \item $\ro(r) \geq 0$ for every $r>c$,
  \item $\lim_{r \to \infty} \ro(r)=\ro< \infty$,
  \item $\lim_{r  \to \infty} r \ro'(r) \log(r) = 0$.
 \end{enumerate}
\end{defi}

\begin{defi}
Two proximate orders $\rho_1(r)$ and $\rho_2(r)$ are said to be \textit{equivalent} if
$$
\lim_{r\to\infty}\big(\rho_1(r)-\rho_2(r)\big)\log(r)=0.
$$
\end{defi}

\begin{rema}\label{remaordenaproxequiv}
If $\rho_1(r)$ and $\rho_2(r)$ are equivalent and $\lim_{r\to\infty}\rho_1(r)=\rho$, then  $\lim_{r\to\infty}\rho_2(r)=\rho$ and  $\lim_{r\to\infty}r^{\rho_1(r)}/r^{\rho_2(r)}=1$.
\end{rema}

From the work of L. S. Maergoiz  we have the following result.

\begin{theo}[\cite{Maergoiz}, Thm.\ 2.4]\label{propanalproxorde}
Let $\ro(r)$ be a proximate order with $\ro(r)\to\ro>0$ as $r\to \infty$. For every $\ga>0$ there exists an analytic function $V(z)$ in $S_\ga$ such that:
  \begin{enumerate}[(i)]
  \item  For every $z \in S_\ga$,
 \begin{equation*}
    \lim_{r \to \infty} \frac{V(zr)}{V(r)}= z^{\ro},
  \end{equation*}
uniformly in the compact sets of $S_\ga$.
\item $\overline{V(z)}=V(\overline{z})$ for every $z \in S_\ga$ (where, for $z=(|z|,\arg(z))$, we put $\overline{z}=(|z|,-\arg(z))$).
\item $V(r)$ is positive in $(0,\infty)$, monotone increasing and $\lim_{r\to 0}V(r)=0$.
\item The function $t\in\R\to V(e^t)$ is strictly convex (i.e. $V$ is strictly convex relative to $\log(r)$).
\item The function $\log(V(r))$ is strictly concave in $(0,\infty)$.
\item  The function $\ro_0(r):=\log( V(r))/\log(r)$, $r>0$, is a proximate order equivalent to $\ro(r)$.
    \end{enumerate}
\end{theo}

We denote by $\mathfrak{B}(\ga,\ro(r))$ the class of such functions $V$. They share a property that will be crucial in the construction of flat functions.

\begin{prop}[\cite{Maergoiz}, Property\ 2.9]\label{propcotaVpartereal}
 Let $\ro>0$, $\ro(r)$ be a proximate order with $\ro(r)\to\ro$, $\ga\ge 2/\ro$ and $V\in \mathfrak{B}(\ga, \ro(r))$. Then, for every $\a\in(0,1/\ro)$ there exist constants $b>0$ and $R_0>0$ such that
 \begin{equation*}
  \Re(V(z)) \ge b V(|z|), \quad  z\in S_{\a},\ |z|\ge R_0,
 \end{equation*}
 where $\Re$ stands for the real part.
\end{prop}

We will also make use of the following result of V. Thilliez.

\begin{theo}[\cite{Thilliez2}, Proposition\ 4]\label{teorcaracfuncplanaAMS}
Let $\M$ be a strongly regular sequence and $S$ a sector. For $f\in\mathcal{H}(S)$, the following are equivalent:
\begin{itemize}
\item[(i)] $f\in\tilde{\mathcal{A}}_{\M}(S)$ and $f\sim_{\M}\hat{0}$.
\item[(ii)] For every bounded proper subsector $T$ of $S$ there exist $c_1,c_2>0$ with
$$|f(z)|\le c_1h_{\M}(c_2|z|)=c_1e^{-M(1/(c_2|z|))},\qquad z\in T.
$$
\end{itemize}
\end{theo}

In the next result we obtain the desired flat functions in case $\omega(\M)<2$ and $d(r)$, defined in Remark~\ref{remadefifuncdder}, is a proximate order. Subsequently, we will indicate how to deal with the case $\omega(\M)\ge 2$. Finally, we will determine conditions on $\M$ amounting to $d(r)$ being a proximate order, or at least guaranteeing that $d(r)$ is a proximate order.

\begin{theo}\label{teorconstrfuncplana}
Suppose $\M$ is a strongly regular sequence with $\omega(\M)<2$ and such that $d(r)$ is a proximate order. Then, for every $V\in\mathfrak{B}(2\omega(\M),d(r))$ the function $G$ defined in $S_{\omega(\M)}$ by
$$
G(z)=\exp(-V(1/z))
$$
belongs to $\tilde{\mathcal{A}}_{\M}(S_{\omega(\M)})$ and it is a (nontrivial) flat function.
\end{theo}

\begin{proof1}
It is enough to reason with sectors $S(0,\omega,r_0)\ll S_{\omega(\M)}$, where $0<\omega<\omega(\M)$ and $r_0>0$. If $z\in S(0,\omega,r_0)$, we have $1/z\in S_{\omega}$.
By our assumptions, $d(r)$ is a proximate order, and by (\ref{equaordeMdet}) and~(\ref{equaordequasMbis}), we have that
$$
\lim_{r\to\infty}d(r)=\ro[M]=\frac{1}{\omega(\M)}.
$$
We are in a position to apply Proposition~\ref{propcotaVpartereal} with $\ro=1/\omega(\M)$,
$\ro(r)=d(r)$, $\ga=2\omega(\M)$ and $\alpha=\omega$, and deduce the existence of constants $R_0>0$ and $b>0$ such that $\Re(V(\zeta)) \ge b V(|\zeta|)$ whenever $\zeta\in S_{\omega}$ with $|\zeta|\ge R_0$. Then, for $z\in S(0,\omega,1/R_0)$ we obtain
$$
|G(z)|=e^{-\Re(V(1/z))}\le e^{-bV(1/|z|)},
$$
and for a suitable $C>0$ we will have $|G(z)|\le Ce^{-bV(1/|z|)}$ for $z\in S(0,\omega,r_0)$.
Now observe that, by the definition of  $\mathfrak{B}(2\omega(\M),d(r))$, we know that the function $\log( V(r))/\log(r)$ is a proximate order equivalent to $d(r)=\log(M(r))/\log(r)$, so that, as a consequence of Remark~\ref{remaordenaproxequiv}, we have that there exists $c>0$ such that for $r> 1/r_0$ one has $V(r)>cM(r)$, and
$$
|G(z)|\le Ce^{-bcM(1/|z|)}=C(h_{\M}(|z|))^{bc}\le Ch_{\M}(D|z|),\qquad z\in S(0,\omega,r_0),
$$
where $D>0$ is a positive constant, suitably chosen according to whether $bc>1$ or not (see  property~(\ref{e120})). It suffices to take into account Theorem~\ref{teorcaracfuncplanaAMS}
 in order to conclude.
\end{proof1}

\begin{rema}\label{remaomegaMmayor2}
In case $\omega(\M)\ge 2$, we may also construct nontrivial flat functions by taking into account the following facts:
\begin{itemize}
\item[(i)] Given a strongly regular sequence $\M=(M_n)_{n\in\N_0}$ and a positive real number $s>0$, the sequence of $s$-powers $\M^{(s)}:=(M_n^s)_{n\in\N_0}$ is strongly regular (see Lemma 1.3.4 in~\cite{thilliez}) and one easily checks that, with self-explaining notation, $\bm^{(s)}=(m_n^s)_{n\in\N_0}$, $M^{(s)}(t)=sM(t^{1/s})$ for every $t\ge 0$, $\omega(\M^{(s)})=s\omega(\M)$, $d^{(s)}(r)=d(r^{1/s})/s+\log(s)/\log(r)$ for $r$ large enough, and
    $$
    r(d^{(s)})'(r)\log(r)=\frac{1}{s}r^{1/s}d'(r^{1/s})\log(r^{1/s})-\frac{\log(s)}{\log(r)}
    $$
    whenever both sides are defined. So, it is clear that $d(r)$ is a proximate order if, and only if, $d^{(s)}(r)$ is.
\item[(ii)] If $\M$ is strongly regular, $\omega(\M)\ge 2$ and $d(r)$ is a proximate order, choose $s>0$ such that $s\omega(\M)<2$. By (i), we may apply Theorem~\ref{teorconstrfuncplana} to $\M^{(s)}$ and obtain $G_0\in\tilde{\mathcal{A}}_{\M^{(s)}}(S_{\omega(\M^{(s)})})$ which is flat. Now, the function $G$, given in $S_{\omega(\M)}$ by $G(z)=G_0(z^s)$, is well-defined and it is plain to see that it is a nontrivial flat element in $\tilde{\mathcal{A}}_{\M}(S_{\omega(\M)})$.
\end{itemize}
\end{rema}

Our next objective is to characterize those $\M$ such that $d(r)$ is a proximate  order.
After looking at Remark~\ref{remadefifuncdder} and Theorems~\ref{teorordeMdetiguallimidder} and~\ref{teororderM}, it is clear that
we only need to care about whether $rd'(r)\log(r)\to 0$ as $r\to\infty$.
The following result provides us with statements equivalent to this fact.

\begin{prop}\label{propcaracdderordenaprox}
Let $\M$ be a strongly regular sequence, and $d(r)$ its associated function. The following are equivalent:
\begin{itemize}
\item[(i)] $d(r)$ is a proximate order,
\item[(ii)] $\lim_{p\to\infty} m_pd'(m_p^+)\log(m_p)=0$,
\item[(iii)] $\displaystyle\lim_{p\to\infty}\frac{p+1}{M(m_{p})}=\frac{1}{\omega(\M)}=\ro[M]$.
\end{itemize}
\end{prop}

\begin{proof1}
For convenience, write $b(r)=rd'(r)\log(r)$ whenever it exists.
It is easy to obtain that
\begin{equation}\label{equaexprrdprimarlogr}
b(r)=\frac{rM'(r)}{M(r)}-d(r)=\frac{p}{M(r)}-d(r),\qquad r\in(m_{p-1},m_p),\ p\in\N.
\end{equation}
Since for sufficiently large $r$, distinct from every $m_p$, we know that $d'(r)>0$, the function $b(r)$ is positive. Moreover, as $M$ and $d$ are both increasing and continuous for large $r$, we see from~(\ref{equaexprrdprimarlogr}) that $b(r)$ is decreasing in every interval $(m_{p-1},m_p)$, and it presents at every $m_p$ a jump of positive height equal to
$$
\lim_{r\to m_p^{+}}b(r)-\lim_{r\to m_p^{-}}b(r)=\frac{1}{M(m_p)}.
$$
From this it is clear that (i) holds if, and only if, $\lim_{p\to\infty}b(m_p^+)=0$, and this is precisely (ii). Now, observe that
$$
b(m_p^+)=\frac{p+1}{M(m_p)}-d(m_p),
$$
and recall from Theorem~\ref{teorordeMdetiguallimidder} that $\lim_{p\to\infty}d(m_p)=\ro[M]=1/\omega(\M)$. So, (ii) amounts to (iii) and we are done.
\end{proof1}

Next we obtain some easy condition that ensures that $d(r)$ is a proximate  order.

\begin{coro}
If
\begin{equation}\label{equacondordenaprox}
\displaystyle\lim_{p\to\infty}p\log\big(\frac{m_{p+1}}{m_{p}}\big) \textrm{ exists (finite or not),}
\end{equation}
then its value is \emph{a fortiori} $\omega(\M)$, $d(r)$ is a proximate order and, moreover, $$\omega(\M)=\lim_{p\to\infty}\frac{\log(m_p)}{\log(p)}\qquad\textrm{ (instead of $\displaystyle\liminf_{p\to\infty}$, see~(\ref{equaordequasM}))}.
$$
\end{coro}

\begin{proof1}
By Stolz's criterion we have that $\lim_{p\to\infty}\frac{\log(m_p)}{\log(p)}$ exists, since
$$
%\lim_{p\to\infty}\frac{\log(m_p)}{\log(p)}=
\lim_{p\to\infty}\frac{\log(m_{p+1}/m_p)}{\log((p+1)/p)}=
\lim_{p\to\infty}p\log\big(\frac{m_{p+1}}{m_p}\big),
$$
and the last limit exists. We take into account~(\ref{equaordequasM}) in order to deduce that
all these limits equal $\omega(\M)$. But, again by Stolz's criterion,
$$
\lim_{p\to\infty}\frac{M(m_p)}{p+1}=\lim_{p\to\infty}\frac{\log(m_p^p/M_p)}{p+1}=
\lim_{p\to\infty}p\log\big(\frac{m_{p}}{m_{p-1}}\big)=\omega(\M),
$$
and this equality amounts to (iii) in Proposition~\ref{propcaracdderordenaprox}.
\end{proof1}

\begin{rema}\label{remacomentdderordenaprox}
\begin{itemize}
\item[(i)] The previous condition~(\ref{equacondordenaprox}) holds for every sequence $\M_{\a,\beta}$, so that in any of these cases $d(r)$ is a proximate order and it is possible to construct flat functions in the corresponding classes. Indeed, we have not been able yet to provide an example of a strongly regular sequence for which $d(r)$ is not a proximate order, i.e., for which condition (iii) in Proposition~\ref{propcaracdderordenaprox} does not hold.
\item[(ii)] In the Gevrey case, $\M_{1/k}=(p!^{1/k})_{\in\N_0}$, let us put $M_{1/k}(r)$, $d_{1/k}(r)$, and so on, to denote the corresponding associated functions. Then, one can check (see, for example, \cite{GelfandShilov}) that for large $r$ we have $c_2r^k\le M_{1/k}(r)\le c_1r^k$ for suitable constants $c_1,c_2>0$, so that $\log(c_2)\le(d_{1/k}(r)-k)\log(r)\le\log(c_1)$ eventually. This shows one can work with the constant proximate order $\ro(r)\equiv k$, and any $V\in\mathfrak{B}(2/k,\ro(r))$ will provide us (due to Theorem~\ref{teorconstrfuncplana}, and since $V(r)$ will be bounded above and below by $r^k$ times some suitable constants) with a flat function in the class $\tilde{\mathcal{A}}_{1/k}(S_{1/k})$. It is easy to see that $V(z)=z^{k}$ belongs to $\mathfrak{B}(2/k,\ro(r))$, and we obtain in this way the classical flat function in this situation, namely $G(z)=\exp(-z^{-k})$.
\item[(iii)] If $\M$ is such that $d(r)$ is not a proximate order, but there exist a proximate order $\ro(r)$ and constants $A,B>0$ such that eventually $A\le (d(r)-\ro(r))\log(r)\le B$, then, by reasoning as indicated in (ii), one may also construct flat functions in $\tilde{\mathcal{A}}_{\M}(S_{\omega(\M)})$.
\end{itemize}
\end{rema}

We are in a position to characterize quasianalyticity in the classes $\tilde{\mathcal{A}}_{\M}(S_{\gamma})$.

\begin{coro}[Watson's Lemma]\label{coroWatsonlemma}
Suppose $\M$ is strongly regular and such that $d(r)$ is a proximate order, and let $\ga>0$ be given. The following statements are equivalent:
\begin{itemize}
\item[(i)] $\tilde{\mathcal{A}}_{\M}(S_{\gamma})$ is quasianalytic, i.e., it does not contain nontrivial flat functions (in other words, the Borel map is injective in this class).
\item[(ii)] $\gamma>\omega(\M)$.
\end{itemize}
\end{coro}

\begin{proof1}
By Theorem \ref{teorconstrfuncplana} and Remark \ref{remaomegaMmayor2}.(ii), whenever $\gamma\le\omega(\M)$ we have nontrivial flat functions in $\tilde{\mathcal{A}}_{\M}(S_{\gamma})$.\par\noindent
Conversely, suppose $\gamma>\omega(\M)$ and that $f\in\tilde{\mathcal{A}}_{\M}(S_{\gamma})$ is a nontrivial flat function. Choose $\gamma'$ with $\omega(\M)<\gamma'<\gamma$. By Proposition \ref{propcotaderidesaasin}.(ii), the restriction of $f$ to $S_{\gamma'}$ belongs to $\mathcal{A}_{\M}(S_{\gamma'})$ and it is flat, so that $\mathcal{A}_{\M}(S_{\gamma'})$ is not quasianalytic, contrary to the definition of $\omega(\M)$.
\end{proof1}

\begin{rema}
One may observe the difference with respect to the classes $\mathcal{A}_{\M}(S_{\gamma})$, which could be quasianalytic for $\gamma=\omega(\M)$ (see Example \ref{examsucesalfabetaomegagamma}).
\end{rema}

\section{Kernels and moment sequences associated with $\bM$}\label{sectkernels}

As a next step in our study, we now devote ourselves to extend to general Carleman classes $\tilde{\mathcal{A}}_{\M}(S_{\gamma})$ the well-known result, named Borel--Ritt--Gevrey theorem, stating that the Borel map in Gevrey classes is surjective if, and only if, the sector is narrow enough. The proof will be constructive, and will rest on the use of truncated Laplace-like transforms whose kernels are intimately related to the nontrivial flat functions obtained in Theorem~\ref{teorconstrfuncplana}. With any such kernel we will associate a sequence of moments which, in turn, will be equivalent to the sequence $\M$ we departed from.

\begin{defi}\label{defikerneleV}
Let $\bM=(M_{p})_{p\in\N_0}$ be a strongly regular sequence such that $d(r)$ is a proximate order, and consider the flat function $G\in\tilde{\mathcal{A}}_{\M}(S_{\omega(\M)})$ constructed in Theorem~\ref{teorconstrfuncplana} for a given $V\in\mathfrak{B}(2\omega(\M),d(r))$. We define the \textit{kernel associated with $V$} as $e_V:S_{\omega(\M)}\to\C$ given by
\begin{equation*}%\label{funcione}
e_V(z):=zG(1/z)=ze^{-V(z)}, \qquad z\in S_{\omega(\M)}.
\end{equation*}
\end{defi}

\begin{rema}\label{remacomentdefinnucleo}
\begin{itemize}
\item[(i)] In a previous paper by A. Lastra, S. Malek and the author~\cite{lastramaleksanz}, similar kernels were obtained from flat functions constructed by V. Thilliez in~\cite{thilliez}. The main difference with respect to the present one, which will be extremely important in forthcoming applications of these ideas to summability theory of formal power series, is that Thilliez needed to slightly restrict the opening of the optimal sector in order to construct such flat functions, while here we have been able to do it in the whole of $S_{\omega(\M)}$.
\item[(ii)] The factor $z$ appearing in $e_V$ takes care of the integrability of $z^{-1}e_V(z)$ at the origin (see $(i)$ in the next lemma). Indeed, it could be changed into any power $z^{\alpha}$ for positive $\a$, where the principal branch of the power is to be considered. Our choice aims at simplicity. %Secondly, as indicated in Remark~\ref{remaConstrG_M}, there are some constants $\delta_1$ and $s$ to be fixed in the construction of $G_{\bM}$.
\end{itemize}
\end{rema}

\begin{lemma}\label{propiedadese}
The function $e_V$ enjoys the following properties:
\begin{itemize}
\item[(i)] $z^{-1}e_V(z)$ is integrable at the origin, it is to say, for any $t_0>0$ and $\tau\in\R$ with $|\tau|<\frac{\pi\omega(\M)}{2}$ the integral $\int_{0}^{t_0}t^{-1}|e_{V}(te^{i\tau})|dt$ is finite.
\item[(ii)] For every $T\prec S_{\omega(\M)}$ there exist $C,K>0$ such that
\begin{equation}\label{e147}|e_V(z)|\le Ch_{\bM}\left(\frac{K}{|z|}\right),\qquad z\in T.
\end{equation}
\item[(iii)] For every $x\in\R$, $x>0$, the value $e_V(x)$ is positive real.
\end{itemize}
\end{lemma}
\begin{proof1}
(i) Let $t_0>0$ and $\tau\in\R$ with $|\tau|<\frac{\pi\omega(\M)}{2}$. Since $G$ is flat,
from Theorem~\ref{teorcaracfuncplanaAMS} we obtain $c_1,c_2>0$ (depending on $\tau$ and $t_0$) such that
$$\int_{0}^{t_0}\frac{|e_V(te^{i\tau})|}{t}dt\le \int_{0}^{t_0}c_1h_{\bM}(c_2/t)dt.$$
As $h_{\bM}$ is continuous and $h_{\bM}(s)\equiv 1$ when $s\ge\frac{1}{m_{1}}$, this integral converges.\par\noindent
(ii) As before, given $T\prec S_{\omega(\M)}$ and $R>0$ there exist $c_1,c_2>0$ (depending on $T$ and $R$) such that
$$|e_V(z)|\le c_1|z|h_{\bM}(c_2/|z|),\qquad z\in T, \ |z|\ge R.$$
If $|z|\ge R$, we may apply (\ref{e120}) for $s=2$ and the definition of $h_{\bM}$ to deduce that
$$
|e_{V}(z)|\le c_1|z|\Big(h_{\bM}\big(\frac{\rho(2)c_2}{|z|}\big)\Big)^2\le c_1|z|h_{\bM}\big(\frac{\rho(2)c_2}{|z|}\big)M_2\big(\frac{\rho(2)c_2}{|z|}\big)^2\le \frac{\rho(2)^2c_1c_2^2M_2}{R}h_{\bM}\big(\frac{\rho(2)c_2}{|z|}\big).
$$
On the other hand, since $V$ is bounded at the origin (because of property (iii) in Theorem~\ref{propanalproxorde}), for $z\in T$ with $|z|<R$ we deduce that
 $e_V(z)=ze^{-V(z)}$ is bounded, and, in order to conclude, it suffices to observe that $h_{\bM}(c_2/|z|)$ is bounded below by some positive constant for $|z|<R$.\par\noindent
(iii) $V(x)$ is real if $x>0$, so $e_V(x)=xe^{-V(x)}>0$.
\end{proof1}

\begin{rema}\label{remanucleoGevreyclasico}
As suggested in Remarks~\ref{remacomentdderordenaprox}.(ii) and~\ref{remacomentdefinnucleo}.(ii), in the Gevrey case $\M_{1/k}$, $k>0$, it is natural and standard to consider the kernel
$$
e_{k}(z)=kz^{k}\exp(-z^{k}),\qquad z\in S_{1/k}.
$$
\end{rema}

\begin{defi}
Let $V\in\mathfrak{B}(2\omega(\M),d(r))$. We define the \textit{moment function associated with $V$} (or to $e_V$) as
$$m_V(\lambda):=\int_{0}^{\infty}t^{\lambda-1}e_{V}(t)dt=
\int_{0}^{\infty}t^{\lambda}G_V(1/t)dt=\int_{0}^{\infty}t^{\lambda}e^{-V(t)}dt.$$
\end{defi}

From Lemma~\ref{propiedadese} we see that $m_V$, well defined in $\{\hbox{Re}(\lambda)\ge0\}$, is continuous in its domain, and holomorphic in $\{\hbox{Re}(\lambda)>0\}$. Moreover, $m_V(x)>0$ for every $x\ge0$. So, the following definition makes sense.

\begin{defi}\label{defisucesmoment}
The sequence of positive real numbers $\mathfrak{m}_V=(m_V(p))_{p\in\N_0}$ is the \textit{sequence of moments associated with $V$} (or to $e_{V}$).
\end{defi}

\begin{prop}\label{mequivm}
Let $e_V$ be a kernel associated with the strongly regular sequence $\M$, and $\mathfrak{m}_V=(m_V(p))_{p\in\N_0}$ the sequence of moments associated with $V$. Then $\M$ and $\mathfrak{m}_V$ are equivalent.
\end{prop}
\begin{proof1}
It suffices to work with $p\ge 1$. From (\ref{e147}) we have $C,K>0$ such that
\begin{equation*}
m_V(p)\le C\int_{0}^{\infty}t^{p-1}h_{\bM}(K/t)dt=
C\int_{0}^{m_{p}}t^{p-1}h_{\bM}(K/t)dt+
C\int_{m_{p}}^{\infty}t^{p-1}h_{\bM}(K/t)dt.
\end{equation*}
In the first integral of the right-hand side we take into account that $h_{\bM}$ is bounded by 1, while in the second one we use the definition of $h_{\bM}$ to obtain that $h_{\bM}(K/t)\le K^{p+1}M_{p+1}/t^{p+1}$, $t>m_p$. This yields
$$m_V(p)\le C\frac{t^p}{p}\Big\vert_{0}^{m_p}- CK^{p+1}M_{p+1}\frac{1}{t}\Big\vert_{m_{p}}^{\infty}=
\frac{Cm_p^p}{p}+CK^{p+1}\frac{M_{p+1}}{m_{p}}.$$
We have $M_{p+1}=m_{p}M_{p}$, and we may apply (\ref{e107}) to obtain that
$$
m_V(p)\le C(A^{2p}+K^{p+1})M_{p}\le 2C\max\{1,K\}(\max\{A^2,K\})^pM_p,
$$
what concludes the first part of the proof.\par
On the other hand, L. S. Maergoiz~\cite[Thm.\ 3.3]{Maergoiz} has shown that the function
\begin{equation}\label{equanucleoE}
F_V(z)=\sum_{n=0}^\infty \frac{z^n}{m_V(n)},\qquad z \in \C,
\end{equation}
is entire and such that
\begin{equation*}
\limsup_{r \to \infty} \frac{\log \max_{|z| = r}|F_V(z)|}{V(r)}\in(0,\infty).
\end{equation*}
From this fact we deduce that there exist constants $C_1,K_1>0$ such that for every $z\in\C$ one has
$$|F_V(z)|\le C_1\exp(K_1V(|z|)).$$
Now, recall that $\log(V(r))/\log(r)$ is a proximate order
equivalent to $d(r)=\log(M(r))/\log(r)$. Consequently, by Remark~\ref{remaordenaproxequiv} there exists $K_2>0$ such that $V(r)\le K_2M(r)$ for large $r$, and so we have
\begin{equation}\label{equadesigF_Vordenaprox}
|F_V(z)|\le \tilde{C}\exp(\tilde{K}M(|z|))
\end{equation}
for every $z\in\C$ and suitably large constants $\tilde{C},\tilde{K}>0$. Finally, we take into account the following result by H. Komatsu~\cite[Prop.\ 4.5]{komatsu}:\par

Let $M(r)$ be the function associated with $\M$. Given an entire function $F(z)=\sum_{n=0}^\infty a_nz^n$, $z\in\C$, the following statements are equivalent:
\begin{itemize}
\item[(i)] There exist $C,K>0$ such that
$|F(z)|\le\displaystyle Ce^{M(K|z|)}$, $z\in\C$.
\item[(ii)] There exist $c,k>0$ such that for every $n\in\N_0$, $|a_n|\le ck^n/M_n$.
\end{itemize}

\noindent
It suffices to apply this equivalence to the function $F_V$, by virtue of~(\ref{equadesigF_Vordenaprox}), and we end the second part of the proof.
\end{proof1}

\begin{rema}
\begin{itemize}
\item[(i)] We record for the future that, as a consequence of the first part of the previous proof, given $K>0$ there exist $C,D>0$ such that for every $p\in\N$ one has
\begin{equation}\label{equacotasintegrhM}
\int_{0}^{\infty}t^{p-1}h_{\bM}(K/t)dt\le CD^pM_p.
\end{equation}
\item[(ii)] In the Gevrey case $\M_{1/k}$ and with the kernel $e_k$ introduced in Remark~\ref{remanucleoGevreyclasico}, we obtain the moment function $m_{1/k}(\lambda)=\Gamma(1+ \lambda/k)$ for $\Re(\lambda)\ge 0$, and we immediately check that $\bM_{1/k}$ and $\mathfrak{m}_{1/k}=(m_{1/k}(p))_{p\in\N_0}$ are equivalent.
\end{itemize}
\end{rema}

\section{A generalization of Borel--Ritt--Gevrey theorem. Right inverses for the asymptotic Borel map}

The proof of the next result, a generalization of the classical Borel--Ritt--Gevrey theorem, will only be sketched, since it is similar to the original one in the Gevrey case (see~\cite{Ramis1,to,ca,balserutx}; in the several variables case, see~\cite{javier}). Indeed, in a previous work by A. Lastra, S. Malek and the author \cite[Thm.\ 4.1]{lastramaleksanz}, this same technique was applied by using kernels derived from the flat functions of V. Thilliez \cite{thilliez}, what obliged us to work in sectors of non-optimal opening. This drawback is now overcome under the additional assumption that the sequence $\bM$ defines a proximate order $d(r)$, which is the case for all the examples we have been able to provide.

\begin{theo}[Generalized Borel--Ritt--Gevrey theorem]\label{tpral}
Let $\bM$ be a strongly regular sequence such that $d(r)$ is a proximate order, and let $\ga>0$ be given. The following statements are equivalent:
\begin{itemize}
\item[(i)] $\gamma\le\omega(\M)$,
\item[(ii)] For every $\boldsymbol{a}=(a_{p})_{p\in\N_0}\in\Lambda_{\bM}$ there exists a function $f\in\tilde{\mathcal{A}}_{\bM}(S_{\gamma})$ such that $$f\sim_{\M}\hat{f}=\sum_{p\in\N_0}\frac{a_p}{p!}z^p,$$
    i.e., $\tilde{\mathcal{B}}(f)=\boldsymbol{a}$. In other words, the Borel map $\tilde{\mathcal{B}}:\tilde{\mathcal{A}}_{\M}(S_{\gamma})\longrightarrow \Lambda_{\M}$ is surjective.
\end{itemize}
\end{theo}

\begin{proof1}
(i)$\implies$(ii) It is enough to treat the case $\ga=\omega(\M)$.
Choose $V\in\mathfrak{B}(2\omega(\M),d(r))$, and consider the associated kernel $e_V$ (see Definition \ref{defikerneleV}) and sequence of moments $\mathfrak{m}_V=(m_V(p))_{p\in\N_0}$ (see Definition \ref{defisucesmoment}).
Given $(a_{p})_{p\in\N_0}\in\Lambda_{\bM}$, there exist $C_1,D_1>0$ such that
\begin{equation*}%\label{equatBoundsLambdaMA}
|a_p|\le C_1D_1^{p}p!M_p,\quad p\in\N_0,
\end{equation*}
so that, by Proposition~\ref{mequivm}, the series
\begin{equation}\label{equatransBorel}
\hat{g}=\sum_{p\in\N_0}\frac{a_p}{p!m_V(p)}z^p
\end{equation}
converges in a disc $D(0,R)$ for some $R>0$, to a holomorphic function $g$. Choose $0<R_0<R$, and define
\begin{equation}\label{intope}
f(z):=\int_{0}^{R_0}e_{V}\left(\frac{u}{z}\right)g(u)\frac{du}{u},\qquad z\in S_{\omega(\M)},
\end{equation}
which turns out to be a holomorphic function in $S_{\omega(\M)}$. Given $T\prec S_{\omega(\M)}$, $N\in\N$ and $z\in T$, by standard arguments we have
\begin{align*}
f(z)-\sum_{p=0}^{N-1}a_p\frac{z^p}{p!} &= f(z)-\sum_{p=0}^{N-1}\frac{a_p}{m_V(p)}m_V(p)\frac{z^p}{p!}\\
&= \int_{0}^{R_0}e_{V}\left(\frac{u}{z}\right)\sum_{k=0}^{\infty}\frac{a_{k}}{m_V(k)}\frac{u^k}{k!}\frac{du}{u} -\sum_{p=0}^{N-1}\frac{a_p}{m_V(p)}\int_{0}^{\infty}u^{p-1}e_{V}(u)du\frac{z^p}{p!}\\
&=\int_{0}^{R_0}e_{V}\left(\frac{u}{z}\right)\sum_{k=N}^{\infty}\frac{a_{k}}{m_V(k)}\frac{u^k}{k!}\frac{du}{u} -\int_{R_0}^{\infty}e_{V}\left(\frac{u}{z}\right)\sum_{p=0}^{N-1}\frac{a_p}{m_V(p)}\frac{u^{p}}{p!}\frac{du}{u}\\
&=f_1(z)+f_2(z).
\end{align*}
By Proposition~\ref{mequivm} there exist $C_2,D_2>0$ such that
\begin{equation}\label{e327}
\frac{|a_{k}|}{m_V(k)k!}\le \frac{C_1D_1^{k}k!M_k}{m_V(k)k!}\le C_2D_2^k
\end{equation}
for all $k\in\N_0$, and so, taking $R_0\le(1-\epsilon)/D_2$ for some $\epsilon>0$ if necessary, we get
\begin{equation}\label{equacotaf1}
|f_{1}(z)|\le C_2\int_{0}^{R_0}\left|e_{V}\left(\frac{u}{z}\right)\right|\sum_{k=N}^{\infty}(D_2u)^{k}\frac{du}{u}
\le \epsilon C_2D_2^{N}\int_{0}^{R_0}\left|e_{V}\left(\frac{u}{z}\right)\right|u^{N-1}du.
\end{equation}
On the other hand, we have $u^p\le R_0^pu^N/R_0^N$ for $u\ge R_0$ and $0\le p\le N-1$. So, according to~(\ref{e327}), we may write
$$\sum_{p=0}^{N-1}\frac{|a_p|u^p}{m_V(p)p!}\le\sum_{p=0}^{N-1}\frac{C_1D_1^pp!M_pu^p}{m_V(p)p!} \le\sum_{p=0}^{N-1}C_1D_1^pC_2D_2^pu^p\le\frac{u^N}{R_0^N}\sum_{p=0}^{N-1}C_1D_1^pC_2D_2^pR_0^p\le C_3D_3^Nu^N$$
for some positive constants $C_3,D_3$, and deduce that
\begin{equation}\label{equacotaf2}
|f_2(z)|\le C_3D^N_3\int_{R_0}^{\infty}\left|e_{V}\left(\frac{u}{z}\right)\right|u^{N-1}du.
\end{equation}
In view of (\ref{equacotaf1}) and (\ref{equacotaf2}), we are done if we prove that
$$
\int_{0}^{\infty}\left|e_{V}\left(\frac{u}{z}\right)\right|u^{N-1}du\le C_4D_4^Nm_V(N)|z|^N
$$
for every $z\in T$ and for suitable $C_4,D_4>0$. But this is a straightforward consequence of Lemma \ref{propiedadese}.(ii) and the estimates in (\ref{equacotasintegrhM}).\par\noindent
(ii)$\implies$(i) We will not provide all the details, but the argument could be completed easily with some of the results in the preprint~\cite{lastramaleksanz2}. Anyway, the idea is similar to the one in the Gevrey case, see~\cite[p.\ 99]{balserutx}. For $\ga>\omega(\M)$, consider a path $\delta_{\omega(\M)}$ in $S_{\ga}$ like the ones used in the classical Borel transform, consisting of a segment from the origin to a point $z_0$ with $\arg(z_0)=\omega(\M)(\pi+\varepsilon)/2$ (for some
$\varepsilon\in(0,\pi)$), then the circular arc $|z|=|z_0|$ from $z_0$ to
the point $z_1$ on the ray $\arg(z)=-\omega(\M)(\pi+\varepsilon)/2$, and
finally the segment from $z_1$ to the origin. Choose any lacunary series $\hat g=\sum_{p=0}^{\infty}b_pz^p/p!$ convergent in the unit disc to a function $g$ that has no analytic continuation beyond that disc (for example, $\hat g=\sum_{p=0}^{\infty}z^{2^p}$). Then, the equivalence of $\M$ and $\mathfrak{m}_V$ implies that $\boldsymbol{a}=(m_V(p)b_p)_{p\in\N_0}$ belongs to $\Lambda_{\M}$. If there would exist a function $f\in\tilde{\mathcal{A}}_{\bM}(S_{\gamma})$ such that $f\sim_{\M}\hat{f}:=\sum_{p\in\N_0}m_V(p)b_pz^p/p!$, one may check that the function
$$
G(u):=\frac{-1}{2\pi i}\int_{\delta_{\omega(\M)}}F_V(u/z)f(z)\frac{dz}{z},\quad
u\in S_{\varepsilon},
$$
where $F_V$ is the function introduced in~(\ref{equanucleoE}), is an analytic continuation of $g$ into the unbounded sector $S_{\varepsilon}$. Since this is not possible, we deduce $\tilde{\mathcal{B}}$ is not surjective in this case.
\end{proof1}

Finally, we will state a result concerning the surjectivity of the asymptotic Borel map $\tilde{\mathcal{B}}$ in the classes $\mathcal{A}_{\M}(S_{\gamma})$, and the existence of suitably defined linear continuous right inverses for $\tilde{\mathcal{B}}$.

\begin{theo}\label{teorinversasderecha}
Let $\M$ be strongly regular and such that $d(r)$ is a proximate order, and let $\ga>0$ be given.\par\noindent
(a) Each of the following assertions implies the next one:
\begin{itemize}
\item[(i)] $\gamma<\omega(\M)$.
\item[(ii)] There exists $d\ge1$ such that for every $A>0$ there is a linear continuous operator
$$
T_{\M,A,\gamma}:\Lambda_{\M,A}\to\mathcal{A}_{\M,dA}(S_{\gamma})
$$
such that $\tilde{\mathcal{B}}\circ T_{\M,A,\gamma}={\emph{Id}}_{\Lambda_{\M,A}}$, the identity map in $\Lambda_{\M,A}$.
\item[(iii)] The Borel map $\tilde{\mathcal{B}}:\mathcal{A}_{\M}(S_{\gamma})\to\Lambda_{\M}$ is surjective.
\item[(iv)] There exists a function $f\in\mathcal{A}_{\M}(S_{\gamma})$ such that for every  $m\in\N_{0}$ we have $f^{(m)}(0)=\delta_{1,m}$ (where $\delta_{1,m}$ stands for Kronecker's delta).
\end{itemize}
\noindent (b) If one has
\begin{equation}\label{equaconddiverenemasunoemeeneomega}
\sum_{n=0}^{\infty}\Big(\frac{1}{(n+1)m_n}\Big)^{1/(\omega(\M)+1)}=\infty,
\end{equation}
then (i) is equivalent to:
\begin{itemize}
\item[(v)] The Borel map $\tilde{\mathcal{B}}:\mathcal{A}_{\M}(S_{\gamma})\to\Lambda_{\M}$ is not injective, i.e., $\mathcal{A}_{\M}(S_{\gamma})$ is not quasianalytic.
\end{itemize}
\noindent (c) If one has
\begin{equation}\label{equaconddiveremeeneomega}
\sum_{n=0}^{\infty}\Big(\frac{1}{m_n}\Big)^{1/\omega(\M)}=\infty,
\end{equation}
then all the conditions (i)-(v) are equivalent to each other.
\end{theo}

\begin{proof1}
(a) (i)$\implies$(ii) Fix $A>0$. For every $\boldsymbol{a}=(a_p)_{p\in\N_0}\in\Lambda_{\M,A}$, the series $\hat g$ given in~(\ref{equatransBorel}) converges in a disc $D(0,R)$ not depending on $\boldsymbol{a}$. We define $T_{\M,A,\gamma}(\boldsymbol{a})$ as the restriction to $S_\ga$ of the function defined in~(\ref{intope}), which was shown to belong to $\tilde{\mathcal{A}}_{\M}(S_{\omega(\M)})$. By combining the information in Proposition~\ref{propcotaderidesaasin} with that in Remark~\ref{remaCarlclassasympexpan},
we conclude that there exists $d\ge1$ such that $T_{\M,A,\gamma}$ sends $\Lambda_{\M,A}$ into $\mathcal{A}_{\M,dA}(S_{\gamma})$ and solves the problem.\par\noindent
(ii)$\implies$(iii) and (iii)$\implies$(iv) are immediate.\par\noindent
(b) By the definition of $\omega(\M)$, we always have that (i) implies (v), and that (v) implies $\ga\le\omega(\M)$. But condition (\ref{equaconddiverenemasunoemeeneomega}) excludes equality by Theorem~\ref{teorKoren}.\par\noindent
(c) Under condition (\ref{equaconddiveremeeneomega}), the fact that (iv) implies (i) may be obtained in the same way as Proposition 3.3 in~\cite{lastrasanz0}. So, (i)-(iv) are all equivalent to each other. According to (b), in order to conclude it suffices to prove that condition (\ref{equaconddiveremeeneomega}) implies condition (\ref{equaconddiverenemasunoemeeneomega}), but this was obtained in Proposition 4.8.(i) in~\cite{lastrasanz1}.
\end{proof1}

\begin{rema}
\begin{itemize}
\item[(i)] Of course, all the results in this paper are valid for general unbounded sectors $S(d,\gamma)$. We have considered the case $d=0$ in the previous arguments only for convenience.
\item[(ii)] V. Thilliez~\cite{thilliez} obtained (i)$\implies$(ii) in the previous result for $\ga<\gamma(\M)$, where $\gamma(\M)$ is the growth index described in Definition~\ref{defiindegrowM}. Since $\gamma(\M)\le\omega(\M)$ in general, our result would mean an improvement for those $\M$ (if any) such that $\gamma(\M)<\omega(\M)$. Also, note that in our present construction of right inverses for $\tilde{\mathcal{B}}$ we need to consider just a ``global'' kernel $e_V$ in $S_{\omega(\M)}$, while in V. Thilliez's and A. Lastra, S. Malek and the author's previous approaches (see~\cite{thilliez,lastramaleksanz}) the kernel had to be chosen depending on the sector $S_\ga$ on which the class was defined.
\item[(iii)] As commented before, the integral expression for the operators $T_{\M,A,\gamma}$ is well suited for their extension to the several variable case. The interested reader may compare this and other approaches in~\cite{javier,lastrasanz,lastramaleksanz}.
\item[(iv)] For Gevrey sequences, condition (\ref{equaconddiveremeeneomega}) holds, since it amounts to the divergence of the harmonic series. In general, condition (\ref{equaconddiverenemasunoemeeneomega}) does not imply  (\ref{equaconddiveremeeneomega}). For instance, as stated in Example~\ref{examsucesalfabetaomegagamma}, the sequence $\M_{\a,\b}$ satisfies (\ref{equaconddiverenemasunoemeeneomega}) if, and only if, $\a\ge\b-1$. One easily checks that it satisfies (\ref{equaconddiveremeeneomega}) if, and only if, $\a\ge\b$. So, if $\b-1\le\a<\b$ we have that $\M_{\a,\b}$ satisfies (\ref{equaconddiverenemasunoemeeneomega}) and not (\ref{equaconddiveremeeneomega}). Whenever this is the case, it is an open problem to decide whether (iv) in the previous theorem implies (i).
\end{itemize}
\end{rema}

\noindent\textbf{Acknowledgements}: This work is partially supported by the Spanish Ministry of Science and Innovation under project MTM2009-12561, and by the Spanish Ministry of Economy and Competitiveness under project MTM2012-31439.

\vskip.5cm
\noindent Author's Affiliation:\par\vskip.5cm

Javier Sanz\par
Departamento de \'Algebra, An\'alisis Matem\'atico, Geometr{\'\i}a y Topolog{\'\i}a\par
Instituto de Investigaci\'on en Matem\'aticas de la Universidad de Valladolid, IMUVA\par
Facultad de Ciencias\par
Universidad de Valladolid\par
47011 Valladolid, Spain\par
E-mail: jsanzg@am.uva.es


\begin{thebibliography}{99}

\bibitem{Balser} W. Balser,
{\it From divergent power series to analytic functions. Theory and application of multisummable power series}, Lecture Notes in Math.
1582, Springer-Verlag, Berlin, 1994.

\bibitem{balserutx} W. Balser, \textit{Formal power series and linear systems of meromorphic ordinary differential equations}, Springer, Berlin, 2000.

\bibitem{bbmt} J. Bonet, R. W. Braun, R. Meise, B. A. Taylor, Whitney's extension theorem for nonquasi-analytic classes of ultradifferentiable functions, Studia Math. {\bf 99} (2) (1991), 155--184.

\bibitem{ca} J. C. Canille, Desenvolvimento assint\'otico e introdu\c{c}\~ao ao c\'alculo diferencial resurgente, 17 Col\'oquio Brasileiro de Matem\'atica, IMPA (1989).

\bibitem{chaucho} J. Chaumat, A. M. Chollet, Surjectivit\'e de l'application restriction \`a un compact dans des classes de fonctions ultradiff\'erentiables, Math. Ann. 298 (1994), no. 1, 7--40.

\bibitem{galindosanz} F. Galindo, J. Sanz, On strongly asymptotically developable functions and the Borel-Ritt theorem, Studia Math. 133 (3) (1999), 231--248.

\bibitem{GelfandShilov} I. M. Gelfand, G. E. Shilov, \textit{Generalized functions, Vol. 2, Space of fundamental and generalized functions}, Academic Press, New York, 1965.

\bibitem{GoldbergOstrowskii} A. A. Goldberg, I. V. Ostrovskii, \textit{Value distribution of meromorphic functions},  Transl. Math. Monogr. 236, Amer. Math. Soc., Providence, RI, 2008.

\bibitem{holland} A. S. B. Holland, \textit{Introduction to the theory of entire functions}, Academic Press, New York and London, 1973.

\bibitem{komatsu} H. Komatsu, Ultradistributions, I: Structure theorems and a characterization, J. Fac. Sci. Univ. Tokyo Sect. IA Math.  20 (1973), 25--105.

\bibitem{korenbljum} B. I. Korenbljum, Conditions of nontriviality of certain classes of functions analytic in a sector, and problems of quasianalyticity, Soviet Math. Dokl. {\bf 7} (1966), 232--236.

\bibitem{lastrasanz0} A. Lastra, J. Sanz, Stieltjes moment problem in general Gelfand-Shilov spaces, Studia Math. 192 (2009), 111--128.

\bibitem{lastrasanz1} A. Lastra, J. Sanz, Quasi-analyticity in Carleman ultraholomorphic classes, Ann. Inst. Fourier 60 (2010), 1629--1648.

\bibitem{lastrasanz} A. Lastra, J. Sanz, Extension operators in Carleman ultraholomorphic classes. J. Math. Anal. Appl. 372 (2010), no. 1, 287--305.

\bibitem{lastramaleksanz} A. Lastra, S. Malek, J. Sanz, Continuous right inverses for the asymptotic Borel map in ultraholomorphic classes via a Laplace-type transform, J. Math. Anal. Appl. 396 (2012), 724--740.

\bibitem{lastramaleksanz2} A. Lastra, S. Malek, J. Sanz, Summability in general Carleman ultraholomorphic classes, submitted. Available at
http://arxiv.org/abs/1402.1669.


\bibitem{Levin} B. Ja. Levin, \textit{Distribution of zeros of entire functions}, Transl. Math. Monogr. 5, Amer. Math. Soc., Providence, RI, 1980.

\bibitem{Maergoiz} L. S. Maergoiz, Indicator diagram and generalized Borel-Laplace transforms for entire functions of a given proximate order, St. Petersburg Math. J. 12 (2001), no. 2, 191--232.

\bibitem{mandelbrojt} S. Mandelbrojt, \textit{S\'eries adh\'erentes, r\'egularisation des suites, applications}, Collection de monographies sur la th\'eorie des fonctions, Gauthier-Villars, Paris, 1952.


\bibitem{MartinetRamis} J. Martinet, J. P. Ramis,
Probl\`{e}mes de modules pour des \'equations diff\'erentielles non
lin\'eaires du premier ordre, Publ. Math. Inst. Hautes Etudes Sci. 55 (1982), 63--164.

\bibitem{pet} H.-J. Petzsche, On E. Borel's theorem, Math. Ann. 282 (1988), no. 2, 299--313.

\bibitem{Ramis1} J. P. Ramis,
D\'evissage Gevrey, Asterisque 59--60 (1978), 173--204.

\bibitem{Ramis2} J. P. Ramis,
{\it Les s\'eries k-sommables et leurs applications}, Lecture Notes in
Phys. 126, Springer-Verlag, Berlin, 1980.

\bibitem{javier} J. Sanz, Linear continuous extension operators for Gevrey classes on polysectors, Glasg. Math. J. 45 (2003), no. 2, 199--216.

\bibitem{SchmetsValdivia} J. Schmets, M. Valdivia, Extension maps in ultradifferentiable and ultraholomorphic function spaces, Studia Math. 143 (3) (2000), 221--250.

\bibitem{Thilliez1} V. Thilliez,
Extension Gevrey et rigidit\'e dans un secteur, Studia Math.
 117 (1995), 29--41.

\bibitem{thilliez} V. Thilliez, Division by flat ultradifferentiable functions and sectorial extensions, Results Math. 44 (2003), 169--188.

\bibitem{Thilliez2} V. Thilliez, Smooth solutions of quasianalytic or ultraholomorphic equations, Monatsh. Math. 160, no. 4 (2010), 443--453.

\bibitem{to} J. Cl. Tougeron, An introduction to the theory of Gevrey expansions and to the Borel-Laplace transform with some applications, Toronto University (1990).

\bibitem{Valiron} G. Valiron, \textit{Th\'eorie des Fonctions}, Masson et Cie., Paris, 1942.

\end{thebibliography}
\end{document}